\documentclass[conference]{IEEEtran}
\IEEEoverridecommandlockouts
\usepackage{graphicx}
\usepackage{amsmath}
\usepackage{amsthm}
\usepackage{amssymb}
\usepackage{epstopdf}
\usepackage{color}
\usepackage{algorithmic,algorithm}
\usepackage{bm}
\usepackage{epstopdf}
\usepackage{adjustbox}
\usepackage{hyperref}

\DeclareMathOperator*{\argmax}{arg\,max}

\newcommand{\R}{\mathbb{R}}

\newcommand{\ndim}{n}
\newcommand{\Ndim}{N}
\newcommand{\nblks}{\ell}

\newcommand{\Ab}{\bm A}
\newcommand{\Wb}{\bm W}

\newcommand{\tb}{\bm t}
\newcommand{\xb}{\bm x}
\newcommand{\yb}{\bm y}

\newcommand{\xbar}{\overline{\xb}}

\newcommand{\XC}{\mathcal X}

\newcommand{\zb}{\bm z}

\newcommand{\essb}{\bm s}

\newcommand{\bb}{\bm b}

\newcommand{\Db}{\bm D}

\newcommand{\Ib}{\bm I}

\newcommand{\proj}{\mathsf{proj}}

\definecolor{darkred}{rgb}{0.75,0,0}

\newcommand{\edit}{\textcolor{blue}}
\renewcommand{\edit}{}

\theoremstyle{plain}
\newtheorem{thm}{Theorem}

\title{A Memory-efficient Algorithm for Large-scale Sparsity Regularized Image Reconstruction}
\author{Greg Ongie, Naveen Murthy, Laura Balzano, Jeffrey A. Fessler
\thanks{
G.~Ongie, N.~Murthy, L.~Balzano, and J.A.~Fessler are with University of Michigan, Dept. of Electrical Engineering and Computer Science, Ann Arbor, MI, USA. E-mail: \{gongie,nnmurthy,girasole,fessler\}@umich.edu.
This work was supported in part by NIH Grant U01 EB018753,  DARPA-16-43-D3M-FP-037, and NSF ECCS-1508943. This work was presented at CT Meeting 2018 (\url{http://www.ct-meeting.org/data/ProceedingsCTMeeting2018.pdf\#page=37}).}}

\begin{document}

\maketitle

\begin{abstract}
 We derive a memory-efficient first-order variable splitting algorithm for convex image reconstruction problems with non-smooth regularization terms. The algorithm is based on a primal-dual approach, where one of the dual variables is updated using a step of the Frank-Wolfe algorithm, rather than the typical proximal point step used in other primal-dual algorithms. We show in certain cases this results in an algorithm with far less memory demand than other first-order methods based on proximal mappings. We demonstrate the algorithm on the problem of sparse-view X-ray computed tomography (CT) reconstruction with \edit{non-smooth} edge-preserving regularization and show competitive run-time with other state-of-the-art algorithms while using much less memory.
\end{abstract}


\section{Introduction}
\label{sec:intro}
Sparsity regularized inverse problems arising in medical imaging result in large-scale non-smooth convex optimization problems that are computationally challenging to solve. 
\edit{General purpose first-order algorithms for nonsmooth convex optimization, such as subgradient descent or smoothing techniques \cite{nesterov:05:smo}, while memory-efficient, converge undesirably slow for this class of problem. Instead, specialized }
proximal splitting algorithms, such as the primal-dual algorithm of Chambolle-Pock (PDCP) \cite{chambolle:11:afo}, represent the current state-of-the-art for these problems. However, PDCP and other proximal methods require storing and operating on one or more dual variables with dimensions potentially several times larger than the image volume to be reconstructed. Even for realistic problem sizes arising in X-ray CT reconstruction, storing and operating on these additional dual variables may be prohibitive or infeasible, limiting the scope of these algorithms in practice. This is especially true when using 3D regularization
with all 26 neighboring voxels
\cite{thibault:07:atd},
or when computing on
GPUs that are limited in memory relative to traditional CPUs. 

To address this issue, this paper introduces a novel algorithm that has far less memory demand than previous approaches.
Specifically, we focus on reconstruction via weighted least squares with \edit{a non-smooth} edge-preserving regularization term. This includes the total variation semi-norm and related penalties. The proposed algorithm is based on a novel primal-dual approach. Existing first-order primal-dual approaches \cite{chambolle:11:afo,condat12aps} alternate between updating the primal variable with a gradient descent step and the dual variable with a projected gradient ascent step (or their proximal equivalents). The main idea of the proposed approach is to replace the projected gradient ascent step in the dual update with a step of the Frank-Wolfe algorithm \cite{frank1956algorithm,jaggi2013revisiting}. We show that this modification allows for substantial memory savings over standard primal-dual approaches. In particular, the algorithm requires storing at most two additional auxiliary variables with dimensions matching the primal variable. We prove convergence of the algorithm under certain assumptions on its step-sizes. 

Finally, we demonstrate the proposed algorithm by reconstructing a sparse-view X-ray CT dataset. Empirically, the proposed algorithm shows competitive convergence with state-of-the-art proximal splitting methods for this problem \cite{sidky:12:cop,nien:15:fxr}, but with much less memory demand.

\vspace{-0.2em}
\section{Problem Formulation}
We consider the following optimization problem:
\vspace{-0.2em}
\begin{equation}\label{eq:primal}
	\min_{\xb} \tfrac{1}{2}\|\Ab\xb - \bb\|_{\Wb}^2 + \lambda\,R(\xb).
\end{equation}
Here $\xb \in \R^n$ represents a vectorized discrete image. The first term in \eqref{eq:primal} measures the data-fit, where $\Ab \in \R^{m\times n}$ with $m \leq n$ is a matrix representing the linear measurement operator, $\bb \in \R^m$ are the (noisy) measurements, and $\|\cdot\|_{\Wb}$ is the weighted $\ell^2$-norm defined as $\|\zb\|_{\Wb} = \sqrt{\zb^T\Wb\zb}$ for a fixed diagonal matrix $\Wb \in \R^{m\times m}$ with positive diagonal entries. The second term is a regularization penalty, where $\lambda > 0$ is a parameter balancing the tightness of data-fit and regularization. This work focuses on regularizers of the form $R(\xb) = \varphi(\Db\xb)$ for some \emph{regularization transform} $\Db\in \R^{\Ndim \times \ndim}$ and where $\varphi(\yb) := \sum_i \phi(\yb_i)$ for some convex and possibly non-smooth sparsity promoting \emph{potential function} $\phi:\R\rightarrow\R^+$. We assume the regularization transform
$\Db$ is a tall matrix ($\Ndim \gg \ndim$) having block form $\Db = [\Db_1^T~\Db_2^T~\cdots~\Db_\nblks^T]^T$ 
with $\Db_i \in \mathbb{R}^{\Ndim_i\times \ndim}$.
For example, if $\{\Db_i\}_{i=1}^\nblks$ is
a collection of first-order finite-difference operators in each dimension and $\phi(x) = |x|$ then $\varphi(\Db\xb) = \|\Db\xb\|_1$ is the discrete (anisotropic) total variation (TV) of $\xb$. Other choices for $\{\Db_i\}_{i=1}^\nblks$ include oriented higher-order finite differences \cite{hu2014generalized} or a collection of pre-trained sparsifying transforms \cite{wen2015structured}. Likewise, the proposed method also generalizes to other convex potential functions $\phi$, such as the Huber loss or Fair potential.

\section{Primal-Dual Frank-Wolfe Algorithm}
For ease of exposition we focus on the case $\varphi(\yb) = \|\yb\|_1$ in the remainder of this work, \emph{i.e.}, $R(\xb) = \|\Db\xb\|_1$. Because the $\ell^\infty$-norm is dual to the $\ell^1$-norm, the primal problem \eqref{eq:primal} with an $\ell^1$ regularizer has an equivalent primal-dual saddle point formulation\footnote{For a general regularization penalty of the form $\varphi(\Db\xb)$, we can derive a similar saddle-point formulation \eqref{eq:primaldual} by writing $\varphi(\Db\xb) = \max_{\yb} \langle \Db\xb,\yb\rangle - \varphi^*(\yb)$ where $\varphi^*$ is the convex conjugate of $\varphi$.} given by
\begin{equation}\label{eq:primaldual}
	\min_{\xb} \max_{\|\yb\|_\infty \leq \lambda} \tfrac{1}{2}\|\Ab\xb - \bb\|_{\Wb}^2 + \langle \Db \xb, \yb\rangle.
\end{equation}
Partitioning the dual vector as
$\yb^T = [\yb_1^T~\cdots~\yb_\nblks^T]$ corresponding to the blocks $\Db^T = [\Db_1^T~\cdots~\Db^T_\nblks]$,
we re-express the inner product in \eqref{eq:primaldual} as
\vspace{-0.5em}
\[
\langle \Db \xb, \yb\rangle = \sum_{i=1}^\nblks \langle \Db_i \xb, \yb_i\rangle =  \bigg\langle \xb,\sum_{i=1}^\nblks\Db_i^T \yb_i\bigg\rangle.
\]
Therefore, introducing the auxiliary variable
$\zb \triangleq \sum_{i=1}^\nblks \Db_i^T\yb_i$, \eqref{eq:primaldual} is equivalent to the equality constrained problem:
\vspace{-0.2em}
\begin{equation}\label{eq:primaldual:3}
	\min_{\xb} \max_{\substack{\|\yb_i\|_\infty \leq \lambda,\\i=1,...,\nblks}} \tfrac{1}{2}\|\Ab\xb - \bb\|_{\Wb}^2 + \left\langle \xb, \zb \right\rangle~~\text{s.t.}~~\zb = \sum_{i=1}^\nblks \Db_i^T\yb_i.
\end{equation}
Below we show that the proposed algorithm only needs to maintain an estimate of the auxiliary variable $\zb \in \R^n$, which has dimensions of the image volume, rather than the full dual variable $\yb \in \R^N$ that is generally several times larger.

We also dualize the data-fit term by defining ${g(\tb) := \frac{1}{2}\|\tb-\bb\|_{\Wb}^2}$ and using the identity
	\[g(\Ab\xb) = \max_{\tb}~\langle\tb,\Ab\xb\rangle - g^*(\tb),\]
where $g^*$ is the convex conjugate of $g$. Simple analysis yields $g^*(\tb) = \frac{1}{2}\|\tb + \Wb\bb\|_{\Wb^{-1}}^2\!\!-\frac{1}{2}\|\bb\|_{\Wb}^2$. Inserting this into \eqref{eq:primaldual:3} and dropping constant terms yields the equivalent formulation
\begin{equation}\label{eq:primaldual:4}
	\min_{\xb} \max_{\tb} \max_{\substack{\|\yb_i\|_\infty \leq \lambda,\\i=1,...,\nblks}} \langle \tb,\Ab\xb\rangle -\frac{1}{2}\|\tb + \Wb\bb\|_{\Wb^{-1}}^2+\langle \xb,\zb \rangle
\end{equation}
~~subject to $\zb = \sum_{i=1}^\nblks \Db_i^T\yb_i$.

\subsection{Frank-Wolfe dual update}
The Frank-Wolfe (FW) algorithm  \cite{frank1956algorithm,jaggi2013revisiting}, also known as the conditional gradient method, is a projection-free approach to solving constrained problems of the form
\[
\max_{\yb\in\mathcal{C}} f(\yb),
\]
where $f$ is a concave function and $\mathcal{C}$ is a closed, convex set. At each iteration, the FW algorithm solves for a search-direction $\essb^\star$ via
\[
\essb^\star = \arg\max_{\essb\in \mathcal{C}}~ \langle \essb, \nabla f(\yb^{(k)})\rangle,
\]
then updates $\yb$ with a convex combination of the previous iterate $\yb^{(k)}$ and the search direction
\[
\yb^{(k+1)} = (1-\alpha_k) \yb^{(k)} + \alpha_k \essb^\star,
\]
where $\alpha_k$ is some iteration-dependent step-size. 

If we apply one step of the FW algorithm to the dual variable $\yb$ in \eqref{eq:primaldual} while holding the primal variable $\xb$ fixed, then the function to maximize is simply the linear function $f(\yb) = \langle \Db \xb, \yb\rangle$ with $\nabla f(\yb) = \Db\xb$ subject to $\|\yb\|_\infty \leq \lambda$. The FW search-direction update in this case is
\[
\essb^\star = \argmax_{\|\essb\|_\infty \leq \lambda }\,\langle\essb, \Db\xb\rangle = \lambda\,\text{sign}(\Db\xb),
\]
\edit{where $\text{sign}(\cdot)$ is applied entrywise and we define $\text{sign}(0)= 0$}.
Hence, a FW update of $\yb$ has the form
\[
\yb^{(k+1)} = (1-\alpha_k) \yb^{(k)} + \alpha_k\,\lambda\,\text{sign}(\Db\xb).
\]
A key to saving memory is that the above update is separable in terms of the $\yb_i$-blocks:
\[
\yb_i^{(k+1)} = (1-\alpha_k) \yb_i^{(k)} + \alpha_k\,\lambda\,\text{sign}(\Db_i\xb), ~\text{for all}~i =1,...,\nblks.
\]
Applying $\Db_i^T$ to both sides above and summing over $i$ yields
\vspace{-1em}
\begin{equation}
\zb^{(k+1)} = (1-\alpha_k) \zb^{(k)} + \alpha_k\lambda\,
\sum_{i=1}^\nblks \Db_i^T\text{sign}(\Db_ix),\label{eq:dualup}
\vspace{-0.5em}
\end{equation}
where we define $\zb^{(k)} = \sum_{i=1}^\nblks \Db_i^T\yb_i^{(k)}$ for all $k\geq 0$.
To save memory, we compute $\zb^{(k+1)}$ incrementally, first by rescaling the current estimate by $(1-\alpha_k)$ then by adding $\alpha_k \lambda \Db_i^T\text{sign}(\Db_i x)$ for all $i=1,...,\nblks$ in sequence.

\subsection{Proximal dual update}
Similar to other primal-dual approaches \cite{chambolle:11:afo,condat12aps}, to update the dual variable $\tb$ we take one step of a proximal point algorithm applied to \eqref{eq:primaldual:4} while fixing the other variables. Specifically, given the current iterates $(\xb^{(k)},\tb^{(k)})$, we set  
\begin{equation*}
\tb^{(k+1)}\!\!=\argmax_{\tb} \langle \tb,\Ab\xb^{(k)} \rangle-\tfrac{1}{2}\|\tb+\Wb\bb\|_{\Wb^{-1}}^2 - \tfrac{1}{2\sigma_k}\|\tb-\tb^{(k)}\|_{\Wb^{-1}}^2
\end{equation*}
where $\sigma_k > 0$ is a step-size parameter to be specified later. This has the closed form solution
\begin{equation}\label{eq:tup}
\tb^{(k+1)} = \tfrac{1}{1+\sigma_k}\tb^{(k)} + \tfrac{\sigma_k}{1+\sigma_k} \Wb(\Ab\xb^{(k)}-\bb).
\end{equation}

\subsection{Primal update}
Finally, we update the primal variable $\xb$ via a gradient descent step (or equivalently a proximal-point step) applied to \eqref{eq:primaldual:4} with the dual variables fixed:
\begin{equation}
\xb^{(k+1)} = \xb^{(k)}-\tau_k(\zb^{(k+1)}+\Ab^T\tb^{(k+1)}),\label{eq:xup}
\end{equation}
where $\tau_k > 0$ is a step-size parameter to be specified later.
Inspired by \cite{chambolle:11:afo} we include an optional over-relaxation step:
\begin{equation}
\xbar^{(k+1)} = \xb^{k+1} + \theta(\xb^{(k+1)}-\xb^{(k)}),\label{eq:xorup}
\end{equation}
where $\theta \in [0,1]$, and perform the dual variable updates \eqref{eq:dualup} and \eqref{eq:tup} with $\xbar^{(k)}$ in place of $\xb^{(k)}$.

\begin{algorithm}
\caption{Primal-Dual Frank-Wolfe (PDFW)}
\label{alg:pdfw}
\begin{algorithmic}
\STATE Initialize $\xb^{(0)} = \xbar^{(0)}, \zb^{(0)} = 0 \in \R^\ndim$ and $\tb^{(0)} = 0 \in \R^{m}$.\\
\STATE Choose step-sizes $(\tau_k,\sigma_k,\alpha_k)$, $\theta\in[0,1]$.
\FOR{$k=0,...,k_{max}$}
	\STATE $\tb^{(k+1)} = \tfrac{1}{1+\sigma_k}\tb^{(k)} + \tfrac{\sigma_k}{1+\sigma_k} \Wb(\Ab\xbar^{(k)}-\bb)$
	\STATE $\zb^{(k+1)} = (1-\alpha_k)\zb^{(k)} + \alpha_k \lambda\sum_{i=1}^\nblks \Db_i^T\text{sign}(\Db_i \xbar^{(k)})$
	\STATE $\xb^{(k+1)} =  \xb^{(k)}\!-\!\tau_k(\Ab^T\tb^{(k+1)} + \zb^{(k+1)})$\\
	\STATE $\xbar^{(k+1)} = \xb^{(k+1)} + \theta(\xb^{(k+1)}-\xb^{(k)})$
\ENDFOR
\end{algorithmic}
\end{algorithm}

\subsection{Algorithm summary and convergence}
Algorithm \ref{alg:pdfw} summarizes the proposed primal-dual Frank-Wolfe (PDFW) algorithm. Using similar analysis as in \cite{bonettini:12:otc} we are able to prove the following convergence result for Algorithm 1 in the special case $\theta = 0$ by showing it is a particular instance of an $\epsilon$-subgradient descent method; we omit the proof for brevity.
\begin{thm}\label{thm:convergence}
Let $\XC^*$ denote the set of minimizers to \eqref{eq:primal}, and let $\{\xb^k\}_{k=1}^\infty$ be the iterates generated by Algorithm 1 with $\theta = 0$. Suppose the iterates $\{\xb^k\}_{k=1}^\infty$ are bounded. If the step-size sequences $\{\alpha_k\} \subset [0,1]$, $\{\sigma_k\} \subset (0,\infty)$, and $\{\tau_k\} \subset (0,\infty)$ satisfy $\tau_k \rightarrow 0,~\sum_{k=0}^\infty \tau_k = + \infty,$ and
\begin{gather*}
	\sum_{j=1}^k \tau_{j-1}\prod_{i=j}^k(1-\alpha_i) \rightarrow 0,~~ 
	\sum_{j=1}^k \tau_{j-1}\prod_{i=j}^k\frac{1}{1+\sigma_k}\rightarrow 0
\end{gather*}
as $k\rightarrow \infty$, then $\mathsf{dist}(\xb^k,\XC^*)\rightarrow 0$, where $\mathsf{dist}$ denotes the Euclidean distance of a point to a set. In particular, if the solution $\xb^*$ to \eqref{eq:primal} is unique then $\xb^k\rightarrow \xb^*$.
\end{thm}
The step-size conditions in Theorem \ref{thm:convergence} are satisfied, for example, when ${\tau_k = O(\frac{1}{k^p})}$, $0 < p \leq 1$, $\alpha_k$ and $\sigma_k$ are constant. There are also valid choices of $\tau_k, \alpha_k,\sigma_k$ for which $\alpha_k \rightarrow 0$ and $\sigma_k \rightarrow \infty$, such as $\tau_k = O(\frac{1}{k^p})$, $\alpha_k = O(\frac{1}{k^q})$, and $\sigma_k = O(\frac{1}{\tau_k})$ with $0 < p < 1$ and $0 < q < p/2$.

Empirically, we observe improved convergence rates using $\theta = 1$ and a constant step-size $\tau_k = \tau$. However, our current proof of Theorem \ref{thm:convergence} does not extend to the case $\theta \neq 0$ nor to the case of $\tau_k$ constant, and we leave its convergence under these conditions as an open problem for future work. 

\subsection{Connections to Chambolle-Pock primal dual algorithm}
Algorithm \ref{alg:pdfw} is closely related to the primal-dual algorithm of Chambolle-Pock (PDCP) \cite{chambolle:11:afo}. If we introduce an auxiliary variable $\essb^{(k)} \in \R^N$ and replace the $\zb^{(k+1)}$ update in Algorithm \ref{alg:pdfw} with the alternative update
\begin{equation*}
	\zb^{(k+1)} = \Db^T\essb^{(k+1)} := \Db^T\proj_{\|\cdot\|_\infty \leq \lambda}(\essb^{(k)}+\sigma_k \Db \xbar^{(k)}),
\end{equation*}
where $\proj_{\|\cdot\|_\infty \leq \lambda}$ denotes Euclidean projection onto the set $\{\essb : \|\essb\|_\infty \leq \lambda \}$, then this modified version of Algorithm 1 coincides with PDCP applied to \eqref{eq:primaldual:4}. In \cite{chambolle:11:afo} it is shown that PDCP converges when $\sigma_k = \sigma$ and $\tau_k = \tau$ are constant and $\tau\sigma L^2 < 1$ and $\theta = 1$, where $L$ is the operator norm of the concatenated matrix $[\Ab^T,\Db^T]^T$.

\subsection{Memory benefits}\label{sec:mem}
Table 1 summarizes the memory requirements of different first-order proximal methods for solving \eqref{eq:primaldual}.
An important feature of Algorithm \ref{alg:pdfw} is that it only requires storing at most three arrays having the size of the image volume to be reconstructed. In contrast, the linearized augmented Lagrangian method (LALM) of \cite{nien:15:fxr} would need to store several arrays have the same size as the image plus two additional arrays of size $N$, the output dimension of the regularization transform. Similarly, the PDCP algorithm \cite{chambolle:11:afo} implemented as in \cite{sidky:12:cop} needs to store at least one array of size $N$. The last column of Table 1 illustrates the memory demand of these algorithms for the iterative reconstruction of a 3D axial CT scan as specified in the next section.
The proposed PDFW algorithm requires an order of magnitude less memory for this example because it avoids having to store large auxiliary variables associated with the regularization transformed image.

\begin{table}
\begin{adjustbox}{width=\columnwidth}
\begin{tabular}{c|c|c|c|c}
& \multicolumn{3}{|c|}{number of variables of size} &  total memory \\
 & $\xb \in \R^n$ & $\Db\xb \in \R^N$ & $\bb \in \R^m$ & 3D CT example\\ 
 & (image) & (reg.~transform) & (data) & (in GB) \\
 \hline
LALM \cite{nien:15:fxr} & 4 & 2 & 2 & 3.02\\ \hline
PDCP \cite{chambolle:11:afo,sidky:12:cop} & 2 & 1 & 2 & 1.60\\ \hline
PDFW, $\theta = 1$ & 3 & 0 & 2 & 0.47 \\ \hline
PDFW, $\theta = 0$ & 2 & 0 & 2 & 0.38 \\ \hline
\end{tabular}
\end{adjustbox}
\vspace{1em}

\caption{Memory demands of first-order methods for solving \eqref{eq:primal}.}
\label{table:mem}
\vspace{-3em}
\end{table}

\begin{figure*}[!ht]
\includegraphics[width=\textwidth]{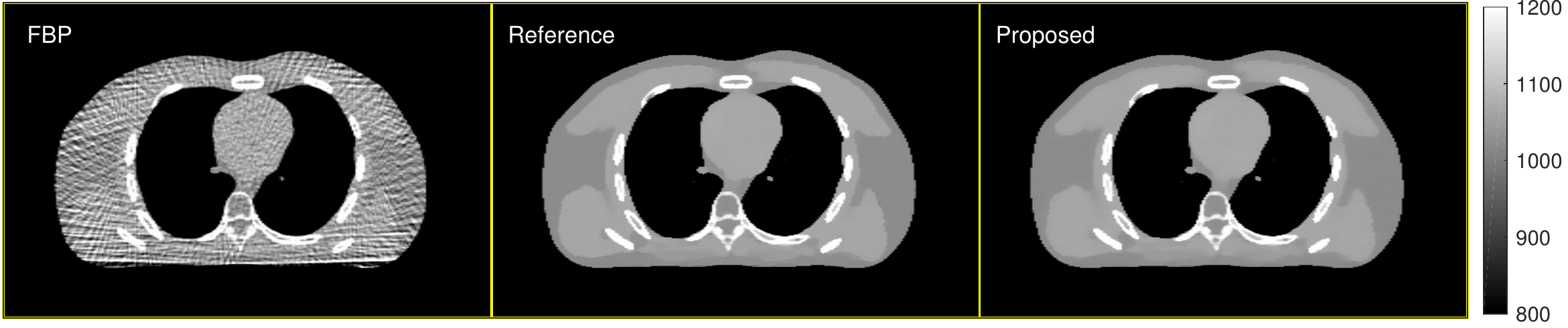}
\caption{Reconstructions of central
transaxial plane of XCAT phantom from sparse-view measurements (left to right): filtered back projection (FBP) reconstruction, a reference solution, and a reconstruction obtained after running 500 iterations of the proposed PDFW algorithm with settings (S2). Images displayed in HU (modified so that air is 0) clipped to range $[800, 1200]$ and cropped to the region of interest.}\label{fig:images}
\end{figure*}
\vspace{-0.2em}

\begin{figure*}[!ht]
\includegraphics[width=\columnwidth]{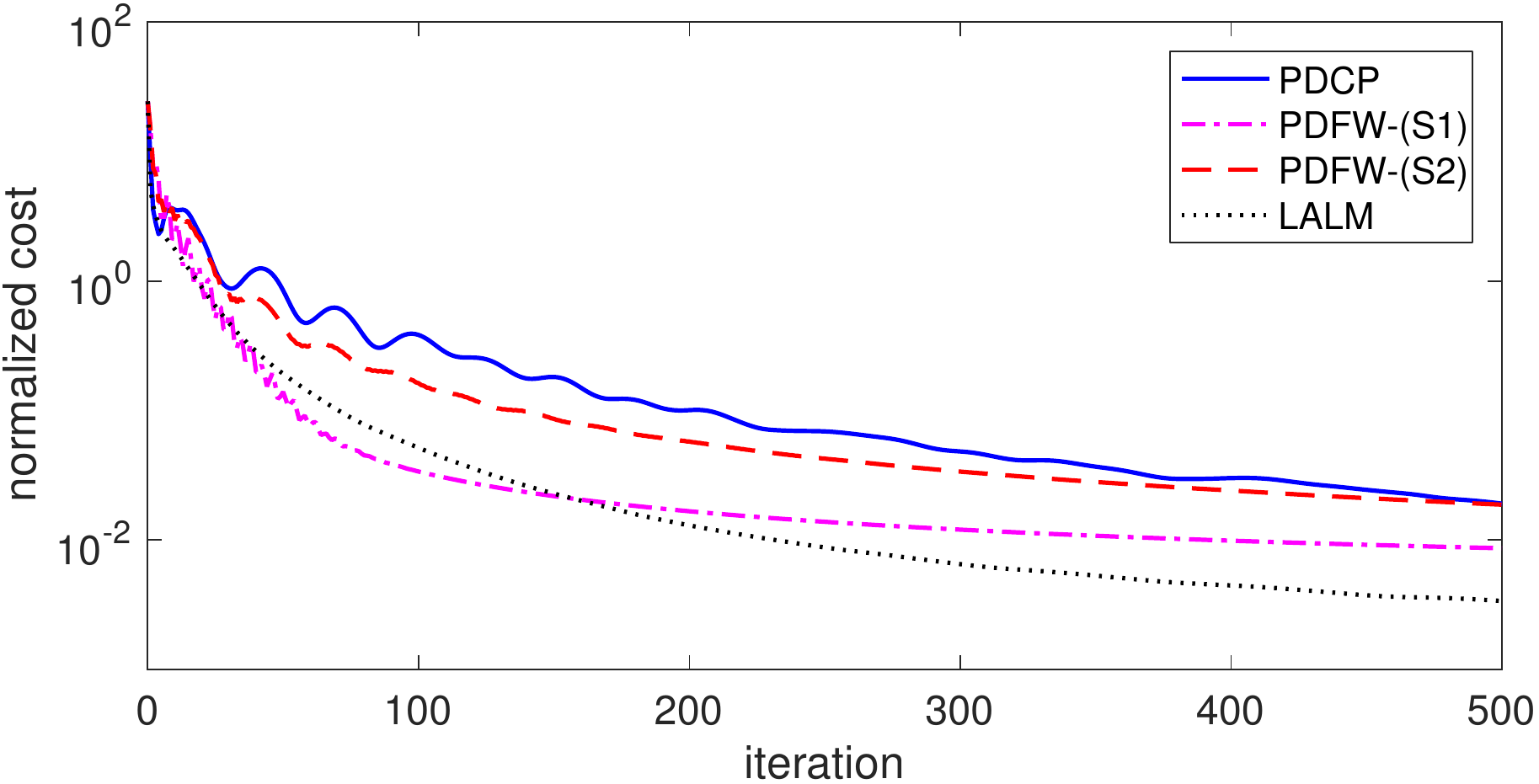}
\includegraphics[width=\columnwidth]{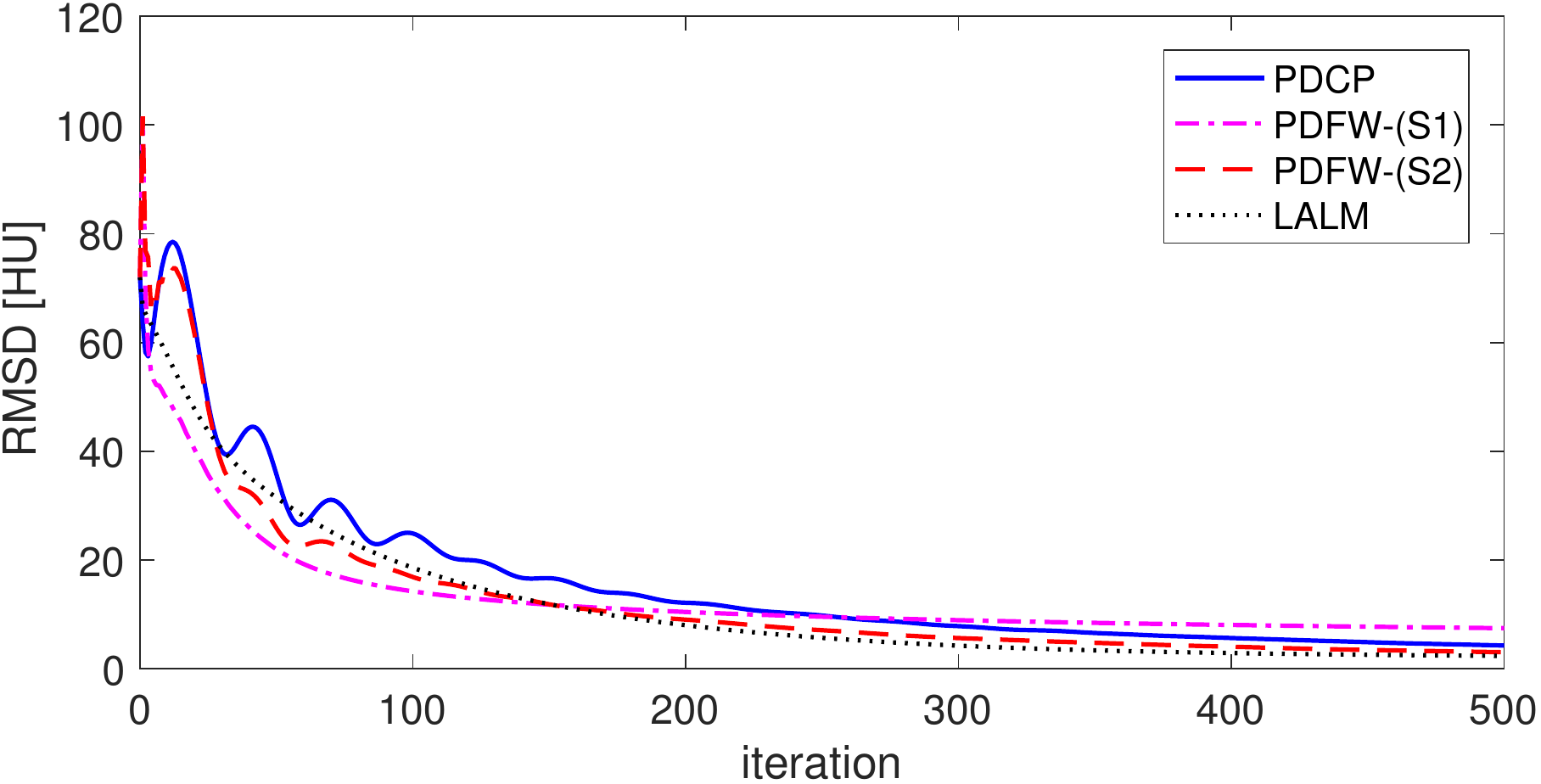}
\caption{Plots of convergence metrics using the proposed PDFW algorithm with settings (S1) and (S2), the PDCP algorithm \cite{sidky:12:cop}, and LALM algorithm \cite{nien:15:fxr}.}
\label{fig:rmsd}
\vspace{-0.7em}
\end{figure*}

\section{Experiments}\label{sec:exp}
Here we demonstrate the proposed PDFW algorithm's potential for sparse-view X-ray CT reconstruction. We simulate an axial CT scan of the XCAT phantom \cite{segars:08:rcs} of size $1024\times 1024\times 154$ voxels to obtain a sinogram of size $m = 888 \times 64\times 120$ (channels $\times$ rows $\times$ views) and reconstruct on a coarser grid of size $n = 512 \times 512 \times 90$. Our reconstruction is obtained by solving \eqref{eq:primal} with regularizer $R(\xb) = \|\Db\xb\|_1$ where $\Db$ computes all finite-differences with thirteen nearest-neighbors of each voxel. We set the statistical weighting matrix $\Wb = \Ib$, and set the regularization parameter $\lambda = 4096$. We compare against two state-of-the-art first-order algorithms for solving \eqref{eq:primal}: PDCP as adapted to CT reconstruction in \cite{sidky:12:cop}, and the linearized augmented Lagrangian method \footnote{An ordered subsets variant of LALM is also presented in \cite{nien:15:fxr}. The proposed PDFW algorithm could also be modified to include ordered subsets updates, but is outside the scope of this work. For fair comparison, we compare against LALM without ordered subsets.} (LALM) of \cite{nien:15:fxr}. 
For the proposed PDFW algorithm we test two sets of step-sizes. The first set (S1) is designed to satisfy the conditions of Theorem 1:
\begin{equation}
	\tau_k = \tfrac{2}{2+k},~~\sigma_k = \tfrac{1}{L^2\tau_k},~~\alpha_k = \left(\tfrac{2}{2+k}\right)^{0.49},~~\theta = 0. \tag{S1}
\end{equation}
The second set (S2) uses a constant step-size $\tau_k=\tau$, violating the conditions of Theorem 1, but matches the settings proposed for the PDCP algorithm in \cite{sidky:12:cop} (except for the choice of $\alpha_k$):
\begin{equation}
	\tau_k = 1/L,~~\sigma_k = 1/L,~~\alpha_k = \tfrac{2}{2+k},~~\theta = 1.\tag{S2}
\end{equation}

Figure \ref{fig:images} shows cropped images from the central
transaxial plane initial filtered back projection reconstruction, reference solution, and the reconstruction obtained from the proposed PDFW algorithm with settings (S2) after 500 iterations. We obtained a reference solution to the optimization problem by running several thousands of iterations of the LALM algorithm, which reached the smallest cost among the competing algorithms. Observe that there is almost no visual difference between the reference solution and the PDFW solution after 500 iterations. 

Figure \ref{fig:rmsd} compares the performance of the algorithms with respect to two convergence metrics: (1) the \emph{normalized cost} defined as ${(f(\xb^{(k)})-f(\xb^*))/f(\xb^*)}$ where $f(\xb)$ is the cost function in \eqref{eq:primal}, $\xb^{(k)}$ is the $k$th iteration of a given algorithm, and $\xb^*$ is the reference solution; and (2) \emph{the root mean square difference} (RMSD) computed as
$\text{RMSD} = \sqrt{\tfrac{1}{|\Omega|}\sum_{{\bm i}\in\Omega}|\xb_{\bm i}^{(k)}-\xb_{\bm i}^*|^2}$
where $\Omega$ is the index set of voxels in a cylindrical region of interest containing the phantom anatomy. 
Overall, the LALM algorithm performs best in terms of the convergence metrics, reaching the lowest cost and RMSD after 500 iterations. However, we reiterate that the LALM algorithm has the highest memory demand of the compared methods (see Table \ref{table:mem}). The proposed PDFW algorithm with step-size scheme (S1) shows a fast initial decrease in the cost and RMSD, but slows in improvement after 100 iterations and has the highest RMSD after 500 iterations, indicating that the (S1) step-size scheme may yield slow asymptotic convergence. The PDFW algorithm with step-size scheme (S2) has better long-run performance in RMSD, yielding nearly the same as LALM after 500 iterations, and its reduction in normalized cost is similar to PDCP. 

\section{Conclusion}
We introduce a memory-efficient algorithm for solving large-scale convex image reconstruction problems with transform sparse regularization based on a novel hybrid of proximal methods and the Frank-Wolfe algorithm. Our experiments demonstrate that the algorithm has competitive performance with other first-order algorithms but with substantially less memory demand. In our experiments we use all of the sinogram measurements to update the primal variable in each iteration, but the proposed algorithm could potentially be modified to incorporate ordered subsets updates similar to \cite{nien:15:fxr} for improved computational efficiency and faster convergence. 

\vspace{-0.3em}
\bibliographystyle{abbrv}

\end{document}